\documentclass[a4paper,12pt]{article}

\usepackage{amssymb,amsmath,latexsym,enumerate,theorem}
\usepackage[all]{xy}

\advance \textheight by \topmargin
\advance \textheight by \headheight
\advance \textheight by \headsep
\advance \textwidth by \marginparsep
\advance \textwidth by \marginparwidth
\advance \textwidth by \oddsidemargin
\def\<{\langle}
\def\>{\rangle}
\def\leq{\leqslant}
\def\geq{\geqslant}

\newcommand{\idR}{{\rm id}_R}

\newcommand{\idP}{{\rm id}_P}
\newcommand{\im}{{\rm im}}
\newcommand{\pr}{{\rm pr}}
\newcommand{\coker}{{\rm coker}}
\newcommand{\Aut}{{\rm Aut}}
\newcommand{\colim}{{\rm colim}}
\newcommand{\GAP}{{\sf GAP}}

\newcommand{\GAPf}{{\sf GAP4}}
\newcommand{\XMod}{{\sf XMod}}
\newcommand{\XModo}{{\sf XMod1}}
\newcommand{\XModt}{{\sf XMod2}}

\newcommand{\calC}{\mathcal{C}}
\newcommand{\calF}{\mathcal{F}}
\newcommand{\calM}{\mathcal{M}}
\newcommand{\sast}{_{\ast}}

\newcommand{\rast}{\; \vec{\ast} \; }

\newcommand{\sbib}{\vspace{-1ex} \bibitem}

\newcommand{\sqdiagram}[8]
  {\xymatrix{ 
    #1 \ar[r]^{#2} \ar[d]_{#4}  &  #3 \ar[d]^{#5} \\
    #6 \ar[r]_{#7}              &  #8 }
  }

\newenvironment{proof}
{\noindent {\bf Proof} }{ \rule{0mm}{2mm} \hfill $\Box$ \mbox{}}

\newtheorem{thm}{Theorem}[section]
{\theorembodyfont{\rmfamily}\newtheorem{example}[thm]{Example}}
\newtheorem{blank}[thm]{\hspace{-0.25em}}
\newtheorem{cor}[thm]{Corollary}
{\theorembodyfont{\rmfamily}\newtheorem{rem}[thm]{Remark}}

\begin{document}

\title{\large\bf 
Computation and Homotopical Applications\\ 
of Induced Crossed Modules}

\author{Ronald Brown \\ Christopher D Wensley \\
{\small Mathematics Division}\\
{\small School of Informatics} \\
{\small University of Wales, Bangor} \\
{\small Gwynedd, LL57 1UT}\\
{\small U.K.}\\
{\small email:~\{r.brown,~c.d.wensley\}@bangor.ac.uk}} 

\maketitle

\begin{center}
{\large Bangor Mathematics Preprint 02.04}
\end{center}
\vspace{1ex}

\begin{abstract}
We explain how the computation of induced crossed modules 
allows the computation of certain homotopy 2-types
and, in particular, second homotopy groups.
We discuss various issues involved in computing induced crossed modules 
and give some examples and applications.
\end{abstract}

\section*{Introduction}

The interactions between topology and combinatorial and computational 
group theory are largely based on the fundamental group functor 
$$
\pi_1 \;:\; \mbox{(based spaces)} \;\to\; \mbox{(groups)}~.
$$ 
At the beginning of the 20th century there was an aim to to generalise
the non commutative fundamental group to higher dimensions, hopes
which seemed to be dashed in 1932 by the proof that the
definition of higher homotopy groups $\pi_n$  then proposed by
\v{C}ech led to commutative groups for $n \geq 2$.

Nonetheless, in the late 1930s and 1940s J.H.C. Whitehead
developed properties of the second relative homotopy group functor
\begin{eqnarray*}
\Pi_2 \;:\; \mbox{(based pairs of spaces)} 
  & \to     &  \mbox{(crossed modules)}~,\\
(X,A,a) 
  & \mapsto & (\partial : \pi_2(X,A,a) \to \pi_1(A,a))~,
\end{eqnarray*}
where $a \in A \subseteq X$  (see Section \ref{topapp}). 
Mac~Lane and Whitehead showed in 1950 \cite{MW} that crossed modules modelled
homotopy 2-types (3-types in their notation) and evidence has
grown that crossed modules can be regarded as `2-dimensional groups'. 
Part of this evidence is the 2-dimensional version of the
Van Kampen Theorem proved by Brown and Higgins in 1978 \cite{BH1},
which allows new computations of homotopy 2-types and so second
homotopy groups. 
This result should be seen as a higher dimensional, non commutative, 
local-to-global theorem,
illustrating themes in Atiyah's  article \cite{atiyah}. 
It is interesting to note that the computation of these second homotopy groups 
is obtained through the computation of a larger non commutative structure. 
This work also throws emphasis on the problem of explicit computation 
with crossed modules, the discussion of which is the theme of this paper.

Our main emphasis in this paper is on induced crossed modules,
which were defined in \cite{BH1} and studied further in papers
by the authors \cite{BW1,BW2}.  
Given the crossed module  $\calM = (\mu: M \to P)$ 
and a morphism of groups  $\iota : P \to Q$, 
the \emph{induced crossed module} $\iota\sast\calM$
has the form $(\partial : \iota\sast M \to Q)$, 
a crossed module over $Q$,
and comes with a morphism of crossed modules
$(\iota\sast, \iota) : \calM \to \iota\sast\calM$~:
$$
\xymatrixcolsep{3pc}  
\xymatrix{
  M \ar[r]^{\iota\sast} \ar[d]_{\mu}
     &  \iota\sast M \ar[d]^{\partial} \\
  P \ar[r]_{\iota} 
     &  ~Q~.
}
$$
Their study requires a solution to many of the general
computational problems of crossed modules. 

In the case $\mu = 0$, when $M$  is simply a $P$-module, 
$\iota\sast M$  is the usual induced $Q$-module 
$M \otimes_{\mathbb{Z}P} \mathbb{Z}Q$. 

Even in the case $M = P,\ \mu = \idP$, 
we know of no relation between the induced crossed module 
$(\partial: \iota_* P\to Q)$ 
and other standard algebraic constructions, although, interestingly,  
$\im\ \partial = N^Q(\iota P)$ the normal closure of  $\iota P$  in  $Q$. 
Thus the induced crossed module construction replaces this
normal closure by a bigger group on which $Q$ acts, and which has
a universal property not usually enjoyed by  $N^Q(\iota P)$.

A long-term project at Bangor is the development of a share library
for the computational group theory program {\GAP} \cite{GAP4}, 
providing functions to compute with these higher-dimensional structures.
The first stage of this project saw the production of the
library {\XModo}, containing functions for crossed modules and
their derivations and for cat$^1$-groups and their sections.
The manual for {\XModo} was included in \cite{GAP3} as Chapter 73.
In particular, Alp \cite{alpthesis} enumerated all isomorphism classes
of cat$^1$-structures on groups of order at most $47$.
This library has recently been rewritten for \GAPf,
with \XModt~ included with the {\sf 4.3} release.
Related libraries include Heyworth's {\sf IdRel} \cite{HW}
for computing identities among the relators of a finitely presented group,
and Moore's {\sf GpdGraph} and {\sf XRes} \cite{moore}
for computing with finite groupoids; group and groupoid graphs;
and crossed resolutions.  These libraries are available at the 
\verb+HDDA+ website \cite{HDDA}.

\section{Crossed modules} \label{xmod}

A \emph{crossed module} $\calM$ (over $P$)
consists of a morphism of groups $\mu : M \to P$,
called the \emph{boundary} of $\calM$,
together with an action of  $P$  on  $M$, written  $(m,p) \mapsto m^p$,
satisfying for all  $m,n \in M, \; p \in P$  the axioms:
$$
CM1)\quad \mu(m^p) = p^{-1}(\mu m) p\;,
\quad\quad
CM2)\quad n^{\mu m} = m^{-1}nm~.
$$
When CM1) is satisfied, but not CM2), 
the structure is a \emph{pre-crossed module} \cite{BHu,HMS}, 
having a \emph{Peiffer subgroup} $C$ 
generated by \emph{Peiffer commutators}
$\< m,n \> \,=\, m^{-1} n^{-1} m\, n^{\mu m}\,$,
and an associated crossed module $(\mu^{\prime} : M/C \to P)$
with  $\mu^{\prime}$  induced by  $\mu$.

Some standard algebraic examples of crossed modules are:
\begin{enumerate}[(i)] 
\item
normal subgroup crossed modules $(i : N \to P)$
where $i$ is an inclusion of a normal subgroup,
and the action is given by conjugation;
\item
automorphism crossed modules $(\chi : M \to \Aut(M) )$
in which $(\chi m)(n) = m^{-1}nm$;
\item
abelian crossed modules $(0 : M \to P)$
where $M$ is a $P$-module;
\item
central extension crossed modules $(\mu : M \to P)$
where $\mu$ is an epimorphism with kernel contained in the centre of $M$.
\end{enumerate}

For our purposes, an important standard construction is the
\emph{free crossed $Q$-module }
$$\calF_{\omega} \;=\; (\partial : F(\omega) \to Q)$$
on a function  $\omega : \Omega \to Q$,
where  $\Omega$  is a set and  $Q$  is a group.
The group  $F(\omega)$ has a presentation with generating set
$\Omega \times Q$  and relators
$$
(m,q)^{-1}\,(n,p)^{-1}\,(m,q)\,(n,pq^{-1}(\omega m)q) \;\;\;
\forall \; m,n \in \Omega, \; p,q \in Q~.
$$
The action is given by  $(m,q)^p = (m,qp)$
and the boundary morphism is defined on generators by
$\partial (m,q) = q^{-1}(\omega m)q$.
This construction will be seen later as a special case of an
\emph{induced} crossed module.
The reader should be warned that the group $F(\omega)$ can be very far
from a free group: in fact, if $\omega$ maps all of $\Omega$ to $\{ 1_Q \}$,
then $F(\omega)$ is just the free $Q$-module on the set $\Omega$,
and in particular is a commutative group.

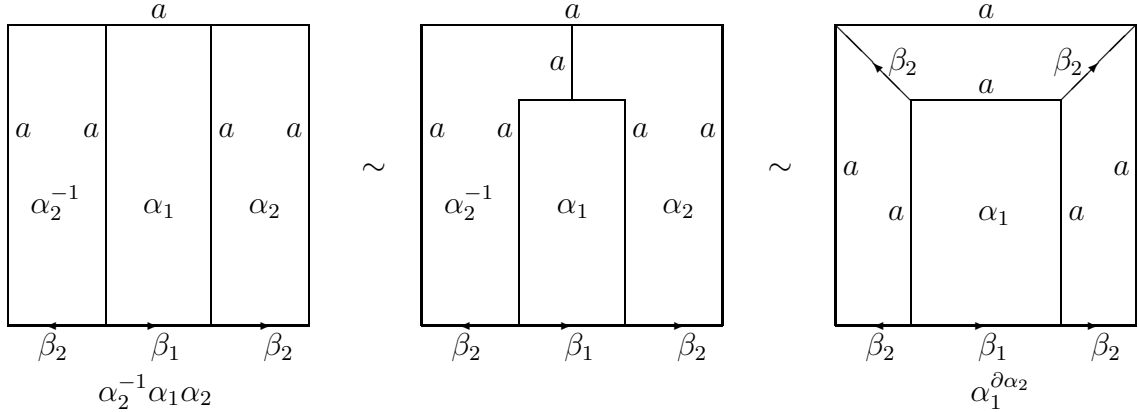
\begin{figure}[!htp] \label{squares}
\setlength{\unitlength}{1mm} 
\begin{center}
\begin{picture}(160,55) 
\put( 5,10){\line(0,1){40}}
\put( 6,35){$a$}
\put( 5,50){\line(1,0){40}}
\put(24,51){$a$}
\put(45,50){\line(0,-1){40}}
\put(42,35){$a$}
\put( 5,10){\line(1,0){5}}
\put(15,10){\vector(-1,0){5}}
\put( 9, 6){$\beta_2$} 
\put(15,10){\vector(1,0){10}}
\put(24, 6){$\beta_1$}
\put(25,10){\vector(1,0){15}}
\put(40,10){\line(1,0){5}}
\put(39, 6){$\beta_2$}
\put(18,10){\line(0,1){40}}
\put(15,35){$a$}
\put(32,10){\line(0,1){40}}
\put(33,35){$a$}
\put(23,25){$\alpha_1$}
\put( 8,25){$\alpha_2^{-1}$}
\put(37,25){$\alpha_2$}
\put( 60,10){\line(0,1){40}}
\put( 61,35){$a$}
\put( 60,50){\line(1,0){40}}
\put( 79,51){$a$}
\put(100,50){\line(0,-1){40}}
\put( 97,35){$a$}
\put( 60,10){\line(1,0){5}}
\put( 70,10){\vector(-1,0){5}}
\put( 64, 6){$\beta_2$} 
\put( 70,10){\vector(1,0){10}}
\put( 79, 6){$\beta_1$}
\put( 80,10){\vector(1,0){15}}
\put( 95,10){\line(1,0){5}}
\put( 94, 6){$\beta_2$}
\put( 73,10){\line(0,1){30}}
\put( 70,35){$a$}
\put( 87,10){\line(0,1){30}}
\put( 88,35){$a$}
\put( 73,40){\line(1,0){14}}
\put( 80,40){\line(0,1){10}}
\put( 77,44){$a$}
\put( 78,25){$\alpha_1$}
\put( 63,25){$\alpha_2^{-1}$}
\put( 92,25){$\alpha_2$}
\put(115,10){\line(0,1){40}}
\put(116,30){$a$}
\put(115,50){\line(1,0){40}}
\put(134,51){$a$}
\put(155,50){\line(0,-1){40}}
\put(152,30){$a$}
\put(115,10){\line(1,0){5}}
\put(125,10){\vector(-1,0){5}}
\put(119, 6){$\beta_2$} 
\put(125,10){\vector(1,0){10}}
\put(134, 6){$\beta_1$}
\put(135,10){\vector(1,0){15}}
\put(150,10){\line(1,0){5}}
\put(149, 6){$\beta_2$}
\put(125,10){\line(0,1){30}}
\put(122,24){$a$}
\put(125,40){\vector(-1,1){5}}
\put(120,45){\line(-1,1){5}}
\put(122,44){$\beta_2$}
\put(125,40){\line(1,0){20}}
\put(134,41){$a$}
\put(145,40){\vector(1,1){5}}
\put(150,45){\line(1,1){5}}
\put(144,44){$\beta_2$}
\put(145,40){\line(0,-1){30}}
\put(146,24){$a$}
\put(134,24){$\alpha_1$}
\put( 52,30){$\sim$}
\put(106,30){$\sim$}
\put( 17,0){$\alpha_2^{-1} \alpha_1 \alpha_2$}
\put(133,0){$\alpha_1^{\partial\alpha_2}$}
\end{picture} 
\caption{\label{secondaxiom} 
~Verification of CM2) for  $\Pi_2(X,A,a)$~.}
\end{center} 
\end{figure}

The major geometric example of a crossed module can be expressed
in two ways.
Let $(X,A,a)$ be a based pair of spaces, with  $a \in A \subseteq X$.
The \emph{second relative homotopy group} $\pi_2(X,A,a)$
consists of homotopy classes rel $J^1$ of continuous maps
$$
\alpha \;:\; (I^2, \dot{I}^2, J^1) \;\to\; (X,A,a)
$$
where  $I=[0,1]$  and  
$J^1=(I \times \{0,1\}) \cup (\{1\} \times I) \subset I^2$.
Each such  $\alpha$  is a map from the unit square  $I^2$  to the space  $X$ 
mapping three sides of the square to the point  $a$  
and the fourth side to a loop at  $a$.
Whitehead showed in \cite{W-46} that there is a crossed module
$\Pi_2(X,A,a)$  with boundary map
$$
\partial \;:\; \pi_2(X,A,a) \to \pi_1(A,a), \quad 
\alpha \mapsto \beta = \alpha(I \times \{0\})~.
$$
The image of  $\alpha_1 \in \pi_2(X,A,a)$
under the action of  $\beta_2 \in \pi_1(A,a)$ 
is illustrated in the right-hand square of Figure \ref{squares}.

Whitehead's main result in \cite{W-41,W-46,W-49} was:
\begin{thm} {\rm (Whitehead)} \label{W} 
If  $X$  is obtained from  $A$  by attaching 2-cells, 
then  $\pi_2(X,A,x)$  is isomorphic to the free crossed
$\pi_1(A,x)$-module on the attaching maps of the 2-cells.
\end{thm}

Later Quillen observed that if $F \to E \to B$ is a based fibration,
then the induced morphism of fundamental groups
$ \pi_1 F \to \pi_1E$ may be given the structure of a crossed module.
This fact is of importance in algebraic $K$-theory.

We also note the following fact, shown in various texts on
homological algebra or the cohomology of groups, e.g. \cite{BK},
and which we relate to topology in section \ref{topapp}:

\begin{blank} \label{k-inv}  
A crossed module  $\calM = (\mu: M \to P)$ 
determines algebraically a cohomology class
$$
k_{\calM} \in H^3(\coker\,\mu, \ker\,\mu),
$$ 
called the $k$-invariant of $\mathcal{M}$, 
and all elements of this cohomology group have such a representation 
by a crossed module.
\end{blank}

\section{Other structures equivalent to crossed modules} \label{other}

One aspect of the problem of higher dimensional group theory is that,
whereas there is essentially only one category of groups,
there are at least five categories of equationally defined algebraic
structures which are equivalent to crossed modules, namely:

\begin{itemize}
\item cat$^1$-groups \cite{L-82};
\item group-groupoids \cite{BS1};
\item simplicial groups with Moore complex of length $1$, \cite{L-82};
\item reduced simplicial $T$-complexes of rank $2$, \cite{dakin,ashley,NT};
\item reduced double groupoids with connection \cite{BS2}.
\end{itemize}

These categories have various geometric models.
The 2-cells of some of these are illustrated in the following pictures:

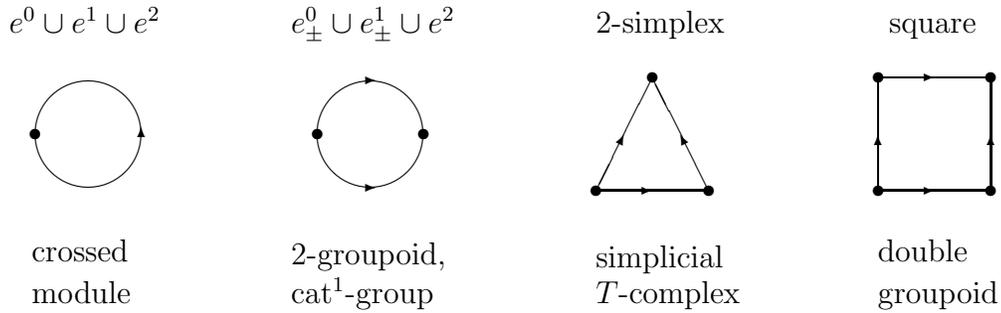
\begin{center}
\setlength{\unitlength}{1.5cm}

\begin{picture}(10,2.5)(-4,-1.5)

\put(-3.2, 0.9)
  {\shortstack[l]{ $e^0 \cup e^1 \cup e^2$ } }
\put(-2.97,0){\circle*{0.1}}
\put(-2.5,0){\circle{1}}
\linethickness{0.7mm}
\put(-2.03,0.07){\vector(0,1){0}}
\put(-3.0,-1.5)
  {\shortstack[l]{ crossed \\
                   module
                   \rule{0mm}{2.5ex} } }

\put(-0.7, 0.9)
  {\shortstack[l]{ $e^0_{\pm} \cup e^1_{\pm} \cup e^2$ } }
\put(-0.47,0){\circle*{0.1}}
\put(0.47,0){\circle*{0.1}}
\put(0,0){\circle{1}}
\linethickness{0.7mm}
\put(0.05,0.475){\vector(1,0){0}}
\put(0.05,-0.475){\vector(1,0){0}}
\put(-0.7,-1.5)
  {\shortstack[l]{ 2-groupoid, \\
                   cat$^1$-group } }

\put(2.0, 0.9)
  {\shortstack[l]{ 2-simplex } }
\put(2,-0.5){\circle*{0.1}}
\put(3,-0.5){\circle*{0.1}}
\put(2.5,0.5){\circle*{0.1}}
\linethickness{0.2mm}
\put(2,-0.5){\line(1,0){1}}
\put(3,-0.5){\line(-1,2){0.5}}
\put(2.5,0.5){\line(-1,-2){0.5}}
\linethickness{0.7mm}
\put(2.5,-0.5){\vector(1,0){0}}
\put(2.75,0.0){\vector(-1,2){0}}
\put(2.25,0){\vector(1,2){0}}
\put(2.0,-1.5)
  {\shortstack[l]{ simplicial \\
                   $T$-complex } }

\put(4.6, 0.9)
  {\shortstack[l]{ square } }
\put(4.5,-0.5){\circle*{0.1}}
\put(5.5,-0.5){\circle*{0.1}}
\put(4.5,0.5){\circle*{0.1}}
\put(5.5,0.5){\circle*{0.1}}
\linethickness{0.2mm}
\put(4.5,-0.5){\line(1,0){1}}
\put(4.5,0.5){\line(1,0){1}}
\put(4.5,-0.5){\line(0,1){1}}
\put(5.5,-0.5){\line(0,1){1}}
\linethickness{0.7mm}
\put(5.0,-0.5){\vector(1,0){0}}
\put(5.0,0.5){\vector(1,0){0}}
\put(4.5,0.0){\vector(0,1){0}}
\put(5.5,0.0){\vector(0,1){0}}
\put(4.5,-1.5)
  {\shortstack[l]{ double \\
                   groupoid
                   \rule{0mm}{2.5ex}  } }
\end{picture}
\end{center}

\noindent
There is also a polyhedral model,
which allows rather general kinds of geometric objects \cite{jones}.

Thus, for computation in ``2-dimensional group theory'',
decisions must be made as to which category to use to represent a
given object, and to compute constructions.
One reason for computing with the crossed module format
is that this is closer to the familiar realm of groups,
for which many computational procedures and systems have
been found and constructed.
Part of the interest in computations with crossed modules is that
such computations will  also yield computations of these
other structures, and this makes them more familiar and understandable.

\subsection{Cat$^1$-groups}

In a cat$^1$-group  $\calC = (e;t,h : G \to R)$
the \emph{embedding}  $e : R \to G$  is a monomorphism
while the \emph{tail} and \emph{head} homomorphisms  $t,h : G \to R$
are surjective and satisfy:
$$
CAT1)\quad  te = he = \idR\,,
\qquad
CAT2)\quad  [\ker t, \ker h] = \{1_G\}\,.
$$
When CAT1) is satisfied, but nor CAT2), 
the structure is a \emph{pre-cat$^1$-group} 
with Peiffer subgroup  $[\ker t, \ker h]$.
A cat$^1$-group  $\calC$  determines a crossed module
$(\partial : S \to R)$  where
$S = \ker t$  and  $\partial = h|_S.$
Conversely, a crossed module $(\mu : M \to P)$
determines a cat$^1$-group
$(e;t,h : P \ltimes M \to P)$ where
$t(p,m) = p$ and $h(p,m) = p(\mu m)$.
The axiom  $he = \idR$ is equivalent to CM1) for a crossed module, 
while CAT2) is equivalent to CM2).
When  $\mu$  is the inclusion of the trivial subgroup in  $P$,
the associated cat$^1$-group  $\calC_P$  has
$e = t = h = \idP$.

Note also that the semidirect product  $P \ltimes M$
admits a groupoid structure with  $t,h$  as source and target,
and composition $\circ$ where
$(p,m) \circ (p(\mu m), n) = (p, mn)$,
making $P \ltimes M$ a group-groupoid,
i.e. a group internal to the category of groupoids.
This notion has a long history: the result that crossed modules
are equivalent to group-groupoids goes back to Verdier,
seems first to have been published in \cite{BS1},
and is used in \cite{breen}.
The holomorph  $\Aut(M) \ltimes M$
of a group $M$ is the source of the cat$^1$-group
associated to the automorphism crossed module
$(\chi : M \to \Aut(M) )$.

Now a colimit  of cat$^1$-groups
$\colim_i(e_i;t_i,h_i : G_i \to R_i)$
is easy to describe.
One takes the colimits $G^{\prime},\ R^{\prime}$
of the underlying groups $G_i,\ R_i$,
and finds that the endomorphisms $e_i,\ t_i,\ h_i$ induce endomorphisms
$e^{\prime} : R^{\prime} \to G^{\prime}$  and 
$t^{\prime},h^{\prime} : G^{\prime} \to R^{\prime}$
satisfying axiom CAT1).
The required colimit is the cat$^1$-group
$\calC^{\prime\prime} =
(e^{\prime\prime}; t^{\prime\prime}, h^{\prime\prime}
   : G^{\prime\prime} \to R^{\prime})$
which has  $G^{\prime\prime} = G^{\prime}/[\ker t^{\prime}, \ker h^{\prime}]$
and  $e^{\prime\prime}, t^{\prime\prime}, h^{\prime\prime}$
induced by  $e^{\prime}, t^{\prime}, h^{\prime}$.

When  $\calC = (e;t,h : G \to R)$  and  $\iota : R \to Q$
is an inclusion, the induced cat$^1$-group  $\iota\sast\calC$
is obtained as the pushout of cat$^1$-morphisms
$(e,\idR) : \calC_R \to \calC$ and
$(\iota,\iota) : \calC_R \to \calC_Q$~.
See Alp \cite{alpthesis}, \cite{AW} for further details.

Further investigation is needed to see whether the use of
cat$^1$-groups can be shown to be more efficient than the direct method 
for the computation of some colimits of crossed modules, 
particularly induced crossed modules,.
The procedure has three stages: convert a crossed module $\calM$
to a cat$^1$-group $\calC$; calculate $\iota\sast\calC$;
then convert $\iota\sast\calC$ to $\iota\sast\calM$.

\section{Computing colimits of crossed modules} \label{colim}

The homotopical reason for interest in computing colimits of
crossed modules is the 2-dimensional Van Kampen Theorem
(2-VKT) due to Brown and Higgins \cite{BH1}. 
The formulation and proof of this theorem was found 
through the notion of double groupoid with connection, 
since such structures yield an appropriate algebraic context 
in which to handle both ``algebraic inverses to subdivision'', 
and the ``homotopy addition lemma''
(which gives a formula for the boundary of a $3$-cube).

One form of the 2-VKT states that Whitehead's \emph{fundamental
crossed module} functor
$$\Pi_2 \;:\; \mbox{(based pairs of spaces)} \to \mbox{(crossed modules)} $$
preserves certain colimits. So for the calculation of certain
homotopy invariants, we need to know how to calculate colimits of
crossed modules. To this end, we start by using some elementary
category theory.

The forgetful functor (crossed modules) $\to$ (groups), 
$(\mu : M \to P) \mapsto P$, has a right adjoint 
$P \mapsto (i : P \to P)$, and so preserves colimits. 
This shows how to compute the $1$-dimensional part of the 
colimit crossed module in terms of colimits of groups. 

The aim now is to transfer the problem to computing colimits of
crossed modules over a fixed group $P$. 
To do this, suppose given a morphism of groups  $\iota : P \to Q $. 
Then there is a pullback functor 
$$
\iota ^{\ast} \;:\; \mbox{(crossed modules over $Q$)} \to
                    \mbox{(crossed modules over $P$)}~.
$$
This functor has a left adjoint 
$$
\iota \sast   \;:\; \mbox{(crossed modules over $P$)} \to  
                    \mbox{(crossed modules over $Q$)}~,
$$ 
which gives our induced crossed module. 
This construction can be described as a ``change of base'' \cite{B-96}.
To compute a colimit  $\colim_i(\mu_i : M_i \to P_i)$, 
one forms the group $P = \mbox{colim}_iP_i$, and
uses the canonical morphisms $\phi _i : P_i \to P$ 
to form the family of induced crossed $P$-modules 
$((\mu_i)\sast : (\phi_i)\sast M_i \to P)$. 
The colimit of these in the category of
crossed $P$-modules is isomorphic to the original colimit. 
Now if  $M^{\prime}$  is the colimit in the category of groups of the
$(\phi_i)\sast M_i$,  then there is a canonical morphism  $M^{\prime} \to P$ 
and an action of $P$ on $M^{\prime}$. 
The resulting  $(M^{\prime} \to P)$  is a pre-crossed module, 
and quotienting by its Peiffer subgroup 
gives the required crossed module.

Presentations for induced crossed modules were given in
\cite{BH1}, and more recently families of explicit examples have
been computed, partly by hand and partly using {\GAP} \cite{BW1}.
Computation of induced crossed modules is here reduced to problems
of computation in combinatorial group theory. 
A key fact which makes one expect successful computations 
is that if $(\mu : M \to P)$ is a crossed module with $M$ finite, 
and if  $\iota : P \to Q$  is a morphism of finite index, 
then the induced crossed $Q$-module
$\iota \sast M$ is also finite \cite[Theorem 2.1]{BW1}.

\begin{example} \label{spec-exs}
When  $\mu : M \to P$  and  $\iota : P \to Q$  are subgroup inclusions, 
there are complete descriptions of  $\iota \sast M$  in the following cases:
\begin{enumerate}[(i)] 
\item  
If  $\iota$  is surjective then $\iota\sast M \cong M /
[M,\ker\iota]$, (\cite[Proposition 9]{BH1}).
\item  
If  $M$  is abelian and  $\iota\mu(M)$  is normal in  $Q$
then  $\iota\sast M$  is abelian and is the usual induced
$Q$-module $M \otimes_{\mathbb{Z}P} \mathbb{Z}Q$, 
(\cite[Corollary 1.6]{BW1}).
\item  
If  $M$  and  $P$  are normal subgroups of  $Q$  then
$\iota\sast M \; \cong \; M \times (M^{{\rm ab}} \otimes I(Q/P))$,
where $I$ denotes the augmentation ideal. 
If in addition  $M=P$  then 
$\iota\sast P \; \cong \; P \times (P^{{\rm ab}})^{[Q:P]-1}$,
(\cite[Theorem 1.1]{BW2}).
\item  
If  $M=P=C_2$, the cyclic group of order $2$,\; $\mu =
\idP$, and  $\iota :C_2 \to D_{2n}$  is the inclusion to a
reflection in the dihedral group $D_{2n}$,  then  $\iota\sast P
\cong D_{2n}$ (\cite[Example 1.4]{BW1}). 
The action is not the usual conjugation:
when  $n$  is odd the boundary is an isomorphism, 
but when  $n$  is even the kernel and cokernel are isomorphic to $C_2$.
\hfill $\Box$
\end{enumerate}
\end{example}

\section{Homotopical applications} \label{topapp}  

As explained in the Introduction, 
the fundamental crossed module functor  $\Pi_2$  assigns  a crossed module 
$(\partial : \pi_2(X,A,a) \to \pi_1(A,a))$ 
to any based pair of spaces $(X,A,a)$. 
Theorem C of \cite{BH1} is a 2-dimensional Van Kampen type theorem 
for this functor.
We will use the following consequence:

\begin{thm} {\em (\cite{BH1}, Theorem D)}  \label{VKT} 
Let $(B,V,b)$ be a cofibred pair of spaces, 
let  $f : V \to A$  be a based map, 
and let  $X$  be the pushout  $A \cup_f B$  in the left-hand diagram below. 
Suppose also that  $A,B,V$  are path-connected, 
and  $(B,V,b)$ is $1$-connected. 
Then the based pair $(X,A,a)$ is $1$-connected and the right-hand diagram 
$$
\xymatrix{
  V \ar[r]^{f} \ar[d]_{\subseteq}
   &  A \ar[d] \\
  B \ar[r]
   &  X
}
\qquad\qquad
\xymatrixcolsep{3pc}  
\xymatrix{
  \pi_2(B,V,b) \ar[r]^{\lambda\sast} \ar[d]_{\delta}
   &  \pi_2(X,A,a) \ar[d]^{\delta^{\prime}} \\
  \pi_1(V,b) \ar[r]_{\lambda} 
   &  \pi_1(A,a)
}
$$
\noindent presents
$\pi_2(X,A,a)$ as the crossed $\pi_1(A,a)$-module
$\lambda_{\ast}(\pi_2(B,V,b))$ induced from the crossed
$\pi_1(V,b)$-module $\pi_2(B,V,b)$ by the group morphism $\lambda :
\pi_1(V,b) \to \pi_1(A,a)$  induced by $f.$  
\end{thm}

As pointed out earlier, when $P$ is a free group on a set $\Omega$  
and $\mu$ is the identity, the induced crossed module
$\iota_{\ast}P$ is the free crossed $Q$-module on the function
$\iota{|_{\Omega}} : \Omega \to Q$. 
Thus Theorem \ref{VKT} implies Whitehead's Theorem 
as stated in Theorem \ref{W}.
A considerable amount of work has been done on this case, 
because of the connections with identities among relations, 
and methods such as transversality theory and ``pictures''  
have proved successful (\cite{BHu,pride}), particularly
in the homotopy theory of 2-dimensional complexes \cite{HMS}.
However, the only route so far available to the wider geometric
applications of induced crossed modules is Theorem \ref{VKT}.  We
also note that this Theorem includes the relative Hurewicz Theorem
in this dimension, on putting $A = \Gamma V$, and $f : V \to
\Gamma V$ the inclusion.  

We will apply this Theorem \ref{VKT} to the {\em classifying space
of a crossed module}, as defined by Loday in \cite{L-82} or Brown
and Higgins in \cite{BH2}.  This classifying space is a functor
$B$ assigning to a crossed module ${\cal M}= (\mu : M \to P)$ a
based $CW$-space $B{\cal M}$  with the following properties:

\noindent   
\begin{blank} 
The homotopy groups of the classifying
space of the  crossed module $\calM = (\mu : M\to P)$ 
are given by
$$
\pi_i(B\calM) \;\;\cong\;\; \left\{  
\begin{array} {rl}
  \coker\, \mu & ~\mbox{\emph{for}~ $i=1$}~, \\ 
    \ker\, \mu & ~\mbox{\emph{for}~ $i=2$}~, \\  
         0     & ~\mbox{\emph{for}~ $i>2$~.}
\end{array} \right. 
$$ 
The first Postnikov invariant of $B{\calM}$
is precisely the $k$-invariant of ${\calM}$ as in \ref{k-inv}.
\end{blank}

\noindent  
\begin{blank}  
The classifying space $BP = B(i : 1 \to P)$ 
is the usual classifying space of the group $P$,  
and $BP$ is a subcomplex of $B\calM$.  
Further, there is a natural isomorphism of crossed modules 
$$
\Pi_2(B\calM, BP, x) \;\;\cong\;\; \calM~.
$$ 
\end{blank}

\noindent  
\begin{blank}  
If $X$ is a reduced $CW$-complex with $1$-skeleton $X^1,$  
then there is a map  
$$
X \;\;\to\;\; B(\Pi_2(X,X^1,x))
$$ 
inducing an isomorphism of $\pi_1$ and $\pi_2$.
\end{blank}

It is in these senses that it is reasonable to say, as in the
Introduction, that crossed modules model all based homotopy 2-types. 

We now give two direct applications of Theorem \ref{VKT}.

\begin{cor} \label{class-sp} 
Let  $\calM = (\mu : M \to P)$ be a crossed module, 
and let $\iota : P\to Q$ be a morphism of groups. 
Let $\beta: BP \to B\calM$ be the inclusion. 
Consider the pushout
$$
\xymatrixcolsep{3pc}  
\xymatrix{
  BP \ar[r]^{\beta} \ar[d]_{B \iota}
     &  B\calM \ar[d] \\
  BQ \ar[r]_{\beta^{\prime}} 
     &  X
}~.
$$
Then the fundamental crossed module of the based pair $(X,BQ,x)$ 
is isomorphic  to the induced crossed module 
$(\partial: \iota_{\ast} M \to Q),$ 
and this crossed module determines the 2-type of $X$. 
In particular, the second homotopy group $\pi_2(X,x)$ 
is isomorphic to  $\ker \partial $.
\end{cor} 
\begin{proof}
The first statement is immediate from Theorem  \ref{VKT}. 
The second statement follows from results of \cite{BH2}, 
since the morphism  $Q \to \pi _1(X)$  is surjective. 
The final statement follows from the homotopy exact sequence of $(X,BQ,x)$.
\end{proof}

\begin{rem} 
An interesting special case of the last
Corollary is when $\calM$ is an inclusion of a normal subgroup, 
since then $B\calM$ is of the homotopy type of  $B(P/M)$.  
So we have determined the 2-type of a homotopy pushout
$$
\sqdiagram{BP}{Bp}{BR}{B\iota}{}{BQ}{p^{\prime}}{X}
$$ 
in which $p : P \to R$ is surjective.
\hfill $\Box$
\end{rem}

\begin{cor}  
Let $\iota :P\to Q$ be a morphism of groups,
and let  $\Gamma BP$  denote the cone on  $BP$.
Then the fundamental crossed module   
$\Pi_2(BQ \cup _{B\iota}\Gamma BP, BQ, x)$  
is isomorphic to the induced crossed module 
$(\partial : \iota \sast P \to Q)$.  
In particular the second homotopy group
$\pi_2 (BQ \cup _{B\iota}\Gamma BP, x)$ 
is isomorphic to  $\ker \partial$. 
\end{cor}

We also note that in determining the crossed module representing a
2-type we are also determining the first Postnikov invariant of that 2-type. 
However it may be more difficult to describe this
invariant as a cohomology class, though this is done in some cases
in \cite{BW1,BW2}.

\section{Computational issues} \label{comp}

We now consider some features of the function {\tt InducedXMod} 
as implemented in {\XMod2}.
The method selection mechanism of {\GAPf} allows for special methods
when  $\iota$  is surjective or injective,
and for the cases listed in example \ref{spec-exs}.

Recall from Proposition 9 of \cite{BH1} that when  $\iota : P \to Q$ 
is a surjection then  $\iota \sast M \cong M/[M,K]$,
where  $K = \ker \iota$  and  $[M,K]$  denotes the subgroup of  $M$ 
generated by the  $m^{-1} m^k$  for all  $m \in M, k \in K$.
When  $\iota$  is neither surjective nor injective, 
we obtain a factorisation  $\iota = \iota_2 \circ \iota_1$ 
with  $\iota_1$  surjective and  $\iota_2$  injective, 
and construct the induced crossed module in two stages:
$$
\xymatrixcolsep{3pc}  
\xymatrix{
  M  \ar[r]^(0.45){(\iota_1) \sast} \ar[d]_{\mu}
  &  (\iota_1) \sast M \ar[r]^(0.55){(\iota_2) \sast} \ar[d]_{\partial_1}
     &  \iota \sast M  \ar[d]^{\partial}        \\
  P  \ar[r]_(0.45){\iota_1} 
  &  \im\,\iota \ar[r]_(0.55){\iota_2} 
     &  ~Q~.
}
$$
The first stage is easily constructed as a quotient group,
so in the following subsections we restrict to the case when 
both  $\iota$  and  $\mu$  are subgroup inclusions.

Note that computation of free crossed modules, 
as described in section \ref{xmod},
is in general difficult since the groups are usually infinite,
and is not attempted in the current version of the package.

\subsection{Copower of groups}

The construction of induced crossed modules, described in \cite{BH1,BW1},
involves the copower  $M \rast T$,
namely the free product  of groups  $M_t, \; t \in T$,
each isomorphic to  $M$.
Here  $T$  is a transversal for the right cosets of  $P$  in  $Q$,
in which the representative of the subgroup  $P$  is taken to be
the identity element.
The group  $M_t = \{ (m,t) \; | \; m \in M \}$
has product  $(m,t)(n,t) = (mn,t)$  and  $Q$  acts by
$(m,t)^q = (m^p,u)$  where  $tq = (\iota p)u$  in  $Q$.
The map  $\delta^{\prime} : M \rast T \to Q$ 
is defined by  $(m,t) \mapsto t^{-1}(\iota \mu m)t$.

The {\GAP} function {\tt IsomorphismFpGroup}
enables the construction of finitely presented groups
$FM, FP, FQ$  isomorphic to the groups  $M, P, Q$; 
monomorphisms  $F\mu : FM \to FP,\; F\iota : FP \to FQ$
mimicing the inclusions  $M \to P \to Q$;
and an action of  $FP$  on  $FM$.
If  $FM$  has  $\gamma$  generators 
then a finitely presented group   $FC$,  isomorphic to  $M \rast T$ 
and with  $\gamma \mid\!T\!\mid$  generators, 
may be constructed using functions in the  $\GAP$
Tietze package (\cite{GAP4}, Chapter 46).
The relators of  $FC$  comprise  $\mid\!T\!\mid$  copies of 
the relators of  $FM$, suitably renumbered.

\subsection{Tracing Tietze transformations}

Let  $\Omega$  be a generating set for  $FM$  and let
$\Omega^{FP}$  be the closure of  $\Omega$  under the action of  $FP$.
Then  $\iota\sast(M) \cong FC/FN$
where  $FN$  is the normal closure in $FC$  of the Peiffer elements
\begin{equation} \label{rels}
  \< (n,s),(m,t) \>
= (n,s)^{-1} (m,t)^{-1} (n,s) (m,t)^{\delta^{\prime}(n,s)} \qquad
  (m,n \in \Omega^{FP}, \; s,t \in T).
\end{equation}
The homomorphism  $\iota\sast$  is induced by
the projection  $\pr_{1}m = (m,1_{FQ})$  onto the first factor,
and the boundary $\delta$  of  $\iota \sast \calM$
is induced from  $\delta^{\prime}$  as shown in the following diagram:
$$
\xymatrixcolsep{3pc}  
\xymatrix{
  FM  \ar[r]^(0.4){\iota\sast} \ar[d]_{\mu}
   &  FC/FN \ar[d]^{\delta} \\
  FP  \ar[r]_(0.4){\iota} 
   &  FQ
}~.
$$
Thus a finitely presented group  $FI \cong \iota \sast M$
is obtained by adding to the relators of  $FC$  further relators
corresponding to the list of elements in equation (\ref{rels}),
and the presentation may be simplified by applying Tietze transformations.

As well as returning an induced crossed module, the construction
should return a morphism of crossed modules
$(\iota \sast, \iota) : \calM \to \iota \sast \calM$.
When Tietze transformations are then applied to the initial presentation
for  $FI$, during the resulting simplification some of the first  $\gamma$  
generators may be eliminated, so the projection  $\pr_{()}$  may be lost.
In order to preserve this projection,
and so obtain the morphism  $\iota \sast$,
it is necessary to record for each eliminated generator  $g$
a relator  $gw^{-1}$  where  $w$  is the word in the remaining generators
by which  $g$  was eliminated.

A significant advantage of {\GAP} is the free availability
of the library code, which enables the user to modify a function
so as to return additional information.
For the {\XModo} version of the package, the Tietze transformation code 
was modified so that the resulting presentation contained an additional field
{\tt presI.remember}, namely a list of (at least)  $\gamma \mid\!T\!\mid$  
relators expressing the original generators in terms of the final ones.
In more recent releases of {\GAP} an equivalent facility has been made
generally available using the {\tt TzInitGeneratorImages} function.

\subsection{Polycyclic groups}

Recall that a polycyclic group is a group  $G$
with power-conjugate presentation having generating set
$ \{ g_1, \ldots, g_n \} $
and relations
\begin{equation}\label{powcon}
   \{ g_i^{o_i} =  w_{ii}(g_{i+1}, \ldots, g_n), \;\;
      g_i^{g_j} = w^{\prime}_{ij}(g_{j+1}, \ldots, g_n)
           \;\;\;  \forall \; 1 \leq j < i \leq n  \}.
\end{equation}
These are implemented in {\GAP} as  {\tt PcGroups} 
(see \cite{GAP4}, Chapters 43,44).
Since subgroups  $M \leq P \leq G$  have induced power-conjugate presentations,
if  $T$  is a transversal for the right cosets of  $P$  in  $G$,
then the relators of  $M \rast T$  are all of the form in (\ref{powcon}).
Furthermore, all the Peiffer relations in equation (\ref{rels})
are of the form  $g_i^{g_j} = g_k^{p}$,
so one might hope that a power conjugate presentation would result.
Consideration of the cyclic-by-cyclic case in the following example
shows that this does not happen in general.

\begin{example}
Let  $C_n$  be cyclic of order  $n$  with generator  $x$,
and let $\alpha : x \mapsto x^a$  
be an automorphism of  $C_n$  of order  $p$. 
Take $G = \< g,h\ |\ g^p, h^n, g^{-1}h^{-1}g\,h^a \> \cong C_p \ltimes C_n$. 
When  $M = P = C_n \lhd G$  cases (ii) and (iii)
of example \ref{spec-exs} apply, and  $\iota\sast C_n \cong C_n^p$. 

It follows from the relators that $h^i g\ =\ g h^{ai},\ 0<i<n,\ $  
and that $h^{-1} (g h^{i(1-a)}) h\ =\ g h^{(i+1)(1-a)}$. 
So if we put  $g_i = gh^{i(1-a)},\ 0 \leq i < n,\ $  
then $g_i^{g_j}\ =\ g_{[j+a(i-j)]}$. 
When  $M = P = C_p = \langle g \ |\ g^p \rangle$  
and  $\iota : C_p \to G$,
we may choose as transversal 
$T = \{ 1_G, h, h^2, \ldots, h^{n-1} \}$.
Then  $M \rast T$  has generators $ \{ (g, h^i)\ |\ 0 \leq i < n \}$, 
all of order $p$, and relators  $ \{ (g,h^i)^p\ |\ 0 \leq i < n \}$. 
The additional Peiffer relators in equation (\ref{rels}) have the form
$$
(g, h^i) (g, h^j) = (g, h^j) (g^k, h^l) 
\quad {\rm where} \quad 
h^i h^{-j} g h^j = g^k h^l
$$
so  $k=1$  and  $l=[j+a(i-j)]$.
Hence  $\theta : \iota\sast M \to Q, \;\; (g,h^i) \mapsto g_i$
is an isomorphism, and  $\iota\sast\calM$
is isomorphic to the identity crossed module on  $Q$.
Furthermore, if we take  $M$  to be a cyclic subgroup  $C_m$  of  $C_p$
then  $\iota\sast\calM$  is the normal subgroup crossed module
$(i : C_m \ltimes C_n \to C_p \ltimes C_n)$.
\hfill $\Box$
\end{example}

\subsection{Identifying  $\iota\sast M$}

From some of the special cases listed in example \ref{spec-exs}
and from other examples, we know that many of the induced groups
$\iota\sast M$ are direct products.  
However the generating sets in the presentations that arise 
following the Tietze transformation do not in general split into 
generating sets for direct summands, 
as the following simple case shows.

\begin{example}
Let $Q=S_4$,  the symmetric group of degree $4$,
and  $M=P=A_4$, the alternating subgroup of $Q$ of index $2$.
Since the abelianisation of  $A_4$  is cyclic of order $3$,
case (iii) in section \ref{colim} shows that
$\iota\sast M \cong A_4 \times C_3$.
However a typical presentation for  $A_4 \times C_3$
obtained from the program is
$$
\< x,y,z\ |\ x^3, y^3, z^3, (xy)^2, zy^{-1}z^{-1}x^{-1},
            yzyx^{-1}z^{-1}, y^{-1}x^2y^2x^{-1}         \>,
$$
and one generator for the $C_3$ summand is  $yzx^2$.

Using the function  {\tt IsomorphismPermGroup}  we obtain 
a permutation group of degree  $12$  with generating set
{\scriptsize
$$
\{
(1,2,3)(4,10,8)(5,11,6)(7,12,9), \;
(1,4,5)(2,6,7)(3,9,10)(8,12,11), \;
(2,5,8)(3,4,7)(6,12,10)
\}.
$$}

On applying  {\tt IsomorphismPcGroup}  to the permutation group 
we obtain a $4$-generator polycyclic group with composition series
$$
A_4 \times C_3 \;\rhd\; A_4 \;\rhd\; C_2^2 \;\rhd\; C_2 \;\rhd\; I~, 
$$
where each subgroup drops the generator  $g_i,\, i=1 \ldots 4$
and  $g_1 g_2 g_4$  is a generator for the normal  $C_3$.
In these representations the cyclic summand  $\ker\partial = C_3$ 
remains hidden, and an explicit search among the normal subgroups 
must be undertaken to find it.
\hfill $\Box$ 
\end{example}

\section{Results}

In this section we list the crossed modules induced from
subgroups of groups of order at most $23$ (excluding $16$), except
that the special cases mentioned earlier enable us to exclude
abelian and dihedral groups; 
the case when  $P$ is normal in $Q$; 
and the case when $Q$ is a semidirect product $C_m \ltimes C_n$.

In the first table, we assume given an inclusion $\iota : P\to Q$
of a subgroup  $P$  of a group  $Q$, and a normal subgroup  $M$  of  $P$. 
We list the isomorphism types of  $\iota \sast M$
and the kernel of  $\partial : \iota_{\ast}M \to Q$. 
Recall that this kernel is realised as a second homotopy group
in corollary \ref{class-sp}.
Labels  $I, C_n, D_{2n}, A_n, S_n$
denote the identity, cyclic, dihedral, alternating and symmetric groups
of order  $1, n, 2n, n!/2$  and  $n!$  respectively.
The group  $H_n$  is the holomorph of  $C_n$
and  $H_n^+$  is its positive subgroup in degree  $n$, 
while  $SL(2,3)$  and  $GL(2,3)$  are the special and general linear groups
of order  $24$  and  $48$.
Labels of the form  {\tt [m,n]}  refer to the $n$th group of order $m$
according to the {\GAPf} numbering.

\begin{center}
\nopagebreak  
\begin{tabular}{||c|c|c|c|c|c||}
\hline
$|Q|$  &  $M$  &  $P$  &  $Q$  &  $\iota\sast M$  &  $\ker \partial$      \\
\hline\hline
$12$ & $C_2$  & $C_2$       & $A_4$              & $H_8^+$     & $C_4$     \\
     & $C_3$  & $C_3$       & $A_4$              & $SL(2,3)$   & $C_2$     \\
$18$ & $C_2$  & $C_2$       & $C_2\ltimes C_3^2$ & $[54,8]$    & $C_3$     \\
     & $S_3$  & $S_3$       & $C_2\ltimes C_3^2$ & $[54,8]$    & $C_3$     \\
$20$ & $C_2$  & $C_2$       & $H_5$              & $D_{10}$    & $C_2$     \\
     & $C_2$  & $C_2^2$     & $D_{20}$           & $D_{10}$    & $I$       \\
     & $C_2^2$& $C_2^2$     & $D_{20}$           & $D_{20}$    & $I$       \\
$21$ & $C_3$  & $C_3$       & $H_7^+$            & $H_7+$      & $I$       \\
\hline
\end{tabular}\\
\vspace{5mm}{\bf Table 1}
\end{center}

The second table contains the results of calculations with  
$Q = S_4$, where  $C_2 = \<(1,2)\>$,  $C_2^{\prime} = \<(1,2)(3,4)\>$, 
and  $C_2^2 = \<(1,2),(3,4)\>$. 
The final column specifiess the automorphism group $\Aut(\iota\sast M)$. 

\begin{center}
\nopagebreak  
\begin{tabular}{||c|c|c|c|c||}
\hline
$M$ & $P$ & $\iota\sast M$ & $\ker \partial$ & $\Aut(\iota\sast M)$  \\
\hline\hline
$C_2$  & $C_2$       & $GL(2,3)$   & $C_2$   & $S_4 C_2$          \\
$C_3$  & $C_3$     & $C_3\ SL(2,3)$& $C_6$   & $[144,183]$        \\
$C_3$  & $S_3$       & $SL(2,3)$   & $C_2$   & $S_4$              \\
$S_3$  & $S_3$       & $GL(2,3)$   & $C_2$   & $S_4 C_2$          \\
$C_2^{\prime}$ & $C_2^{\prime}$  & $[128,?]$ & $C_4C_2^3$  &      \\
$C_2^{\prime}$ & $C_2^2,C_4$ & $H_8^+$     & $C_4$   & $S_4 C_2$  \\
$C_2^{\prime}$ & $D_8$       & $C_2^3$     & $C_2$   & $SL(3,2)$  \\
$C_2^2$& $C_2^2$     & $S_4C_2$    & $C_2$   & $S_4 C_2$          \\
$C_2^2$& $D_8$       & $S_4$       & $I$     & $S_4$              \\
$C_4$  & $C_4$       & $[96,219]$    & $C_4$   & $[96,227]$        \\
$C_4$  & $D_8$       & $S_4$       & $I$     & $S_4$              \\
$D_8$  & $D_8$       & $S_4 C_2$   & $C_2$   & $S_4 C_2$          \\
\hline
\end{tabular}\\
\vspace{5mm}{\bf Table 2}
\end{center}

An interesting problem is to obtain a clearer
understanding of the geometric significance of these tables.


{\small

\begin{thebibliography}{99}

\sbib{alpthesis} {\sc Alp, M.},
\emph{GAP, Crossed modules, Cat$^1$-groups:
     Applications of computational group theory}
Ph.D. thesis, University of Wales, Bangor, (1997),\\
\verb+http://www.informatics.bangor.ac.uk/public/math/research/ftp/theses/alp.ps.gz+.

\sbib{AW} {\sc Alp, M. and Wensley, C.D.},
`Enumeration of cat$^1$-groups of low order',
\emph{Int. J. Algebra and Computation} 10 (2000) 407-424.

\sbib{ashley} {\sc Ashley, N.K.},
\emph{Simplicial $T$-complexes},
Ph.D. Thesis, University of Wales, Bangor, (1978).
Published as
`Simplicial T-complexes: a non-abelian version of a theorem of Dold-Kan',
\emph{Diss. Math.} 265 (1988) 11-58.

\sbib{atiyah} {\sc Atiyah, M.}, 
`Mathematics in the 20th Century',
\emph{Bull. London Math. Soc.} 34 (2002) 1-15.

\sbib{breen} {\sc Breen, L.},
`Th\'eorie de Schreier sup\'erieure',
\emph{Ann. Sci. \'Ecol. Norm. Sup.} 25 (1992) 465-514.

\sbib{BK} {\sc Brown, K.S.},  
\emph{Cohomology of groups},
Graduate texts in Mathematics 87, 
Springer-Verlag, New York (1982).

\sbib{B-96} {\sc Brown, R.}, 
`Homotopy theory, and change of base for groupoids and multiple groupoids', 
\emph{Applied categorical structures}, 4 (1996) 175-193.

\sbib{BH1} {\sc Brown, R. and Higgins, P.J.}, 
`On the connection between the second relative homotopy groups
 of some related spaces', 
\emph{Proc. London Math. Soc.}, (3) 36 (1978) 193-212.

\sbib{BH2} {\sc Brown, R. and Higgins, P.J.},  
`The classifying space of a crossed complex',  
{\em Math. Proc. Camb. Phil. Soc.} 110 (1991) 95-120.

\sbib{BHu} {\sc Brown, R. and Huebschmann, J.},
`Identities among relations',
in  \emph{Low-dimensional topology},
ed. {\sc R.Brown and T.L.Thickstun},
London Math. Soc. Lecture Note Series 46,
Cambridge University Press (1982) 153-202.

\sbib{BS1} {\sc Brown, R. and Spencer, C.B.},
`Double groupoids and crossed modules',
\emph{Cah.  Top. G\'eom. Diff.}, 17 (1976) 343-362.

\sbib{BS2} {\sc Brown, R. and Spencer, C.B.},
`$\mathcal{G}$-groupoids, crossed modules
and the fundamental groupoid of a topological group',
\emph{Proc. Kon. Ned. Akad. v. Wet.}, 79 (1976)  296-302.

\sbib{BW1}  {\sc Brown, R. and Wensley, C.D.},
`On finite induced crossed modules
and the homotopy 2-type of mapping cones',
\emph{Theory and applications of categories}, 1 (1995) 51-74.

\sbib{BW2}  {\sc Brown, R. and Wensley, C.D.},
`Computing crossed modules induced by an inclusion of a normal subgroup,
with applications to homotopy theory',
\emph{Theory and applications of categories}, 2 (1996) 3-16.

\sbib{dakin} {\sc K. Dakin},
\emph{Multiple compositions for higher dimensional groupoids},
Ph.D. Thesis, University of Wales, Bangor, (1977).

\sbib{GAP4}
{\sc The GAP~Group}, 
\emph{{GAP -- Groups, Algorithms, and Programming}}, 
Version 4.3 (2002), \verb+http://www.gap-system.org+.

\sbib{HW}
{\sc Heyworth, A. and Wensley, C.D.},
`Logged rewriting and identities among relators',
in  \emph{Proceedings of Groups St Andrews 2001 in Oxford},
(to appear).

\sbib{HDDA}
{\sc Higher-dimensional discrete algebra},\\
\verb+http://www.informatics.bangor.ac.uk/public/math/research/hdda/+.

\sbib{HMS}
{\sc Hog-Angeloni, C., Metzler, W. and Sieradski, A.J. (Editors)},
\emph{Two-dimensional homotopy and combinatorial group theory},
London Math. Soc. Lecture Note Series 197,
Cambridge University Press (1993).

\sbib{jones}  {\sc Jones, D.W.},
\emph{Polyhedral $T$-complexes},
Ph.D. thesis, University of Wales, Bangor, (1984).
Published as  
`A general theory of polyhedral sets and their corresponding $T$-complexes',
\emph{Diss. Math.} 266 (1988).

\sbib{L-82} {\sc Loday, J.-L.},
`Spaces with finitely many non-trivial homotopy groups',
\emph{J. Pure Appl. Algebra}, 24 (1982) 179-202.

\sbib{MW} {\sc Mac Lane, S. and Whitehead, J.H.C.}, 
`On the 3-type of a complex',  
{\em Proc. Nat. Acad. Sci.} (1950) 41-48.

\sbib{moore} {\sc Moore, E.J.},
\emph{Graphs of groups: word computations and free crossed resolutions}, 
Ph.D. thesis, University of Wales, Bangor, (2001),\\
\verb+http://www.informatics.bangor.ac.uk/public/math/research/ftp/theses/moore.ps.gz+.

\sbib{NT} {\sc Nan Tie, G.},
`A Dold-Kan theorem for crossed complexes',
\emph{J. Pure Appl. Algebra}, 56 (1989) 177-194.

\sbib{pride} {\sc Pride, S.J.}
`Identities among relations',
in  \emph{Proc. Workshop on Group Theory from a Geometrical Viewpoint,
International Centre of Theoretical Physics, Trieste, 1990},
ed. {\sc E. Ghys, A. Haefliger and A. Verjodsky},
World Scientific (1991) 687-716.

\bibitem{GAP3} {\sc { Sch{\accent127 o}nert, M. et~al}},
\emph{GAP--Groups, Algorithms, and Programming},
Lehrstuhl D f{\accent127 u}r Mathematik,
Rheinisch Westf{\accent127 a}lische Technische Hoch\-schule,
Aachen, Germany, sixth edition, 1997.

\sbib{W-41} {\sc Whitehead, J.H.C.},
`On adding relations to homotopy groups',
\emph{Annals of Math.} 41 (1941) 806-810.

\sbib{W-46} {\sc Whitehead, J.H.C.},
`Note on a previous paper entitled
``On adding relations to homotopy groups''',
\emph{Annals of Math.} 47 (1946) 806-810.

\sbib{W-49} {\sc Whitehead, J.H.C.},
`Combinatorial homotopy II',
\emph{Bull. Amer. Math. Soc.} 55 (1949) 453-496.

\end{thebibliography}

}
\end{document}